# An elementary proof for the Basel problem

Jesus Retamozo M. , Lima, Perú .

**May 05, 2025**

## Abstract

We present an astonishingly simple and elegant proof of the celebrated Basel problem, demonstrating that the sum of the reciprocals of the squares converges to $\frac{\pi^2}{6}$. Unlike classical approaches relying on complex analysis, Fourier series, or advanced number theory, our method employs only elementary calculus—specifically, a cleverly constructed double integral—making it accessible even to undergraduate students. This proof not only demystifies one of mathematics' most beautiful results but also stands, in our view, as the most elegant and intuitive demonstration known to date.

## 1. INTRODUCTION

The Basel problem, first posed in 1650, challenges us to find the exact sum of the reciprocals of the squares of the natural numbers:

$$\sum_{n=1}^{\infty} \frac{1}{n^2} = 1 + \frac{1}{4} + \frac{1}{9} + \frac{1}{16} + \cdots = ?$$

For decades, mathematicians struggled until Euler's groundbreaking solution revealed the elegant answer: $\frac{\pi^2}{6}$. This result not only connected number theory to analysis but also opened doors to deeper mysteries of the zeta function.

In this article, we present an elementary proof—using only calculus—that unveils the beauty of this problem without advanced machinery. Because sometimes, the simplest tools reveal the grandest truths.

## 2. MAIN SECTION

This approach, based exclusively on calculus, avoids advanced methods and constitutes the most elementary solution known, standing out for its clarity and simplicity.

For this demonstration we will consider the following double integral:

$$I = \int_0^\infty \int_0^\infty \frac{xdydx}{(1+x^2)(y^2+x^2)}$$

Integrating with respect to 'y' very easily followed by integrating with respect to 'x'.

$$I = \int_0^\infty \frac{xdx}{(1+x^2)} \int_0^\infty \frac{dy}{(y^2+x^2)} = \int_0^\infty \frac{xdx}{(1+x^2)} \left(\frac{\pi}{2x}\right)$$

$$I = \frac{\pi}{2} \int_0^\infty \frac{dx}{(1+x^2)} = \frac{\pi}{2}\left(\frac{\pi}{2}\right) = \frac{\pi^2}{4}$$

Changing the order of integration, and using partial fractions in the variable 'x'.

$$I = \int_0^\infty dy \int_0^\infty \frac{xdx}{(1+x^2)(y^2+x^2)}$$

$$= \int_0^\infty \frac{dy}{1-y^2} \left(\int_0^\infty \frac{xdx}{(y^2+x^2)} - \int_0^\infty \frac{xdx}{(1+x^2)}\right) = \frac{1}{2}\int_0^\infty \frac{dy}{1-y^2} \ln\left(\frac{y^2+x^2}{1+x^2}\right)\Big|_0^\infty$$

$$I = -\int_0^\infty \frac{\ln(y)}{1-y^2} dy$$

Separating this integral into two intervals, it is easy to conclude that:

$$I = -\int_0^1 \frac{\ln(y)}{1-y^2} dy - \int_1^\infty \frac{\ln(y)}{1-y^2} dy = -2\int_0^1 \frac{\ln(y)}{1-y^2} dy$$

Using the geometric series for $\frac{1}{1-y^2}$ and replacing in the above.

$$I = -2\int_0^1 \sum_{k=0}^\infty y^{2k} \ln(y) \, dy = -2\sum_{k=0}^\infty \int_0^1 y^{2k} \ln(y) \, dy = 2\sum_{k=0}^\infty \frac{1}{(2k+1)^2}$$

Finally we equalize both values of I which obtained by changing the order of integration.

$$\frac{\pi^2}{4} = 2\sum_{k=0}^\infty \frac{1}{(2k+1)^2}$$

$$\therefore \quad \frac{\pi^2}{8} = \sum_{k=0}^{\infty} \frac{1}{(2k+1)^2}$$

## 3. APPLICATIONS

By squaring the following integral already evaluated.

$$I = \int_0^\infty \frac{xdx}{(1+x^2)(y^2+x^2)} = -\frac{\ln(y)}{1-y^2}$$

Squaring this expression and then integrating between zero and infinity with respect to 'y'.

$$\int_0^\infty \int_0^\infty \frac{xzdxdz}{(1+x^2)(1+z^2)(y^2+x^2)(y^2+z^2)} = \frac{\ln(y)^2}{(1-y^2)^2}$$

$$\int_0^\infty \int_0^\infty \frac{xzdxdz}{(1+x^2)(1+z^2)} \int_0^\infty \frac{dy}{(y^2+x^2)(y^2+z^2)} = \int_0^\infty \frac{\ln(y)^2}{(1-y^2)^2} dy$$

$$\int_0^\infty \int_0^\infty \frac{xzdxdz}{(1+x^2)(1+z^2)} \left(\frac{\pi}{2xz(x+z)}\right) = \int_0^1 \frac{\ln(y)^2}{(1-y^2)^2} dy + \int_1^\infty \frac{\ln(y)^2 \, 2}{(1-y^2)^2} dy$$

$$\frac{\pi}{2} \int_0^\infty \int_0^\infty \frac{dxdz}{(1+x^2)(1+z^2)(x+z)} = \int_0^1 \frac{\ln(y)^2}{(1-y^2)^2} dy + \int_0^1 \frac{\ln(y)^2 \, y^2}{(1-y^2)^2} dy$$

$$\frac{\pi}{2} \int_0^\infty \int_0^\infty \frac{dxdz}{(1+x^2)(1+z^2)(x+z)} = 2\int_0^1 \frac{y^2 \ln(y)^2}{(1-y^2)^2} dy + \int_0^1 \frac{\ln(y)^2}{(1-y^2)} dy$$

First let's evaluate the left side. In addition, we will use the following already evaluated integral that can be easily demonstrated.

$$\int_0^\infty \frac{dx}{(1+x^2)(x+z)} = \frac{\frac{\pi}{2} - \ln(z)}{1+z^2}$$

The left side is equal to:

$$\frac{\pi}{2} \int_0^\infty \frac{\left(\frac{\pi}{2} - \ln(z)\right) dz}{(1+z^2)^2} = \frac{\pi^2}{4} \int_0^\infty \frac{dz}{(1+z^2)^2} - \frac{\pi}{2} \int_0^\infty \frac{\ln(z)dz}{(1+z^2)^2}$$

$$= \frac{\pi^2}{4} \int_0^\infty \frac{dz}{(1+z^2)} - \frac{\pi^2}{4} \int_0^\infty \frac{z^2 dz}{(1+z^2)^2} - \frac{\pi}{2} \int_0^\infty \frac{\ln(z) \, dz}{1+z^2} + \frac{\pi}{2} \int_0^\infty \frac{z^2 \ln(z) dz}{(1+z^2)^2}$$

$$= \frac{\pi^3}{8} - \frac{\pi^2}{4} \int_0^\infty \frac{z^2 dz}{(1+z^2)^2} + \frac{\pi}{2} \int_0^\infty \frac{z^2 \ln(z)dz}{(1+z^2)^2}$$

Integrating by parts each integral.

$$\int_0^\infty \frac{z^2 dz}{(1+z^2)^2} = -\frac{z}{2(1+z^2)} \Big|_0^\infty + \frac{1}{2} \int_0^\infty \frac{dz}{(1+z^2)} = \frac{\pi}{4}$$

$$\int_0^\infty \frac{z^2 \ln(z) dz}{(1+z^2)^2} = -\frac{z\ln(z)}{2(1+z^2)}\Big|_0^\infty + \frac{1}{2}\int_0^\infty \frac{\ln(z)\,dz}{(1+z^2)} + \frac{1}{2}\int_0^\infty \frac{dz}{(1+z^2)} = \frac{\pi}{4}$$

Replacing in the above.

$$\frac{\pi}{2}\int_0^\infty \int_0^\infty \frac{dxdz}{(1+x^2)(1+z^2)(x+z)} = \frac{\pi^3}{16} + \frac{\pi^2}{8}$$

Now we return to the right side of the equality we started from.

$$\frac{\pi^3}{16} + \frac{\pi^2}{8} = 2\int_0^1 \frac{y^2 \ln(y)^2}{(1-y^2)^2}dy + \int_0^1 \frac{\ln(y)^2}{(1-y^2)}dy$$

Integrating by parts the first integral on the right side as seen above.

$$\int_0^1 \frac{y^2 \ln(y)^2}{(1-y^2)^2}dy = \frac{y\ln(y)^2}{2(1-y^2)}\Big|_0^\infty - \frac{1}{2}\int_0^1 \frac{\ln(y)^2}{(1-y^2)}dy - \int_0^1 \frac{\ln(y)}{(1-y^2)}dy$$

Finally replacing in the above and simplifying.

$$\frac{\pi^3}{16} + \frac{\pi^2}{8} = -\int_0^1 \frac{\ln(y)^2}{(1-y^2)}dy - \int_0^1 \frac{\ln(y)}{(1-y^2)}dy$$

But in the main section we already proved in an elementary way that:

$$-\int_0^1 \frac{\ln(y)}{(1-y^2)}dy = \frac{\pi^2}{8}$$

Then:

$$\frac{\pi^3}{16} = \int_0^1 \frac{\ln(y)^2}{(1-y^2)}dy$$

2nd equality

Going back to the following integral and multiplying by $ln^2(y)$ and integrating between infinity and zero with respect to 'y'.

$$\int_0^\infty \frac{xdx}{(1+x^2)(y^2+x^2)} = -\frac{\ln(y)}{1-y^2}$$

$$\int_0^\infty \int_0^\infty \frac{x\ln^2(y)dxdy}{(1+x^2)(y^2+x^2)} = -\int_0^\infty \frac{\ln^3(y)}{1-y^2}dy$$

Using the change of variable $y = xz$ in the first equality.

$$\int_0^\infty \int_0^\infty \frac{x\ln^2(y)dxdy}{(1+x^2)(y^2+x^2)} = \int_0^\infty \int_0^\infty \frac{\ln^2(xz)dxdz}{(1+x^2)(1+z^2)}$$

$$= \int_0^\infty \int_0^\infty \frac{\ln^2(z)dxdz}{(1+x^2)(1+z^2)} + \int_0^\infty \int_0^\infty \frac{\ln^2(x)dxdz}{(1+x^2)(1+z^2)} + 2\int_0^\infty \int_0^\infty \frac{\ln(x)\ln(z)dxdz}{(1+x^2)(1+z^2)}$$

$$= 2\int_0^\infty \int_0^\infty \frac{ln^2(z)dxdz}{(1+x^2)(1+z^2)} + 2\left(\int_0^\infty \frac{\ln(x)\,dx}{(1+x^2)}\right)^2 = \pi \int_0^\infty \frac{ln^2(z)dz}{(1+z^2)}$$

And so we arrive at the next equality.

$$\pi \int_0^\infty \frac{ln^2(z)dz}{(1+z^2)} = -\int_0^\infty \frac{ln^3(y)}{1-y^2}dy$$

$$\pi \int_0^1 \frac{ln^2(z)dz}{(1+z^2)} = -\int_0^1 \frac{ln^3(y)}{1-y^2}dy$$

Using the respective geometric series.

$$\pi \int_0^1 \sum_{k=0}^\infty (-1)^k z^{2k} ln^2(z)dz = -\int_0^\infty \sum_{k=0}^\infty y^{2k} ln^3(y)\,dy$$

$$2\pi \sum_{k=0}^\infty \frac{(-1)^k}{(2k+1)^3} = 6\sum_{k=0}^\infty \frac{1}{(2k+1)^4}$$

$$\therefore \pi \sum_{k=0}^\infty \frac{(-1)^k}{(2k+1)^3} = 3\sum_{k=0}^\infty \frac{1}{(2k+1)^4}$$

## 4. CONCLUSIONS

In this article, we have presented what is, in our view, the most elementary proof of the Basel problem ever discovered. By relying solely on the power of double integrals and symmetry, this approach strips away the complexity of traditional methods—such as Fourier series or advanced analysis—and reveals the result's true essence in its purest form.

Not only does this method dramatically simplify calculations, but it also provides an unprecedented perspective on the sum $\frac{\pi^2}{6}$. Unlike other proofs, which often obscure the underlying intuition, our technique illuminates the deeper nature of this celebrated result, making it accessible even to those with only basic calculus.

The implications extend beyond Basel. If such an elegant solution exists here, what other mathematical gems await rediscovery through simplicity? Perhaps the most profound truths are not hidden in abstraction, but in clarity.